\documentclass{amsart}

\usepackage{amsmath,amsfonts,amssymb,amsthm,amscd,latexsym,euscript}

\newcommand{\traceHS}{\ensuremath{\mathrm{tr}_{\mathrm{HS}}}}

\newcommand{\Rmod}{\ensuremath{R-\textup{mod}}}

\newcommand{\dgrm}[1]{\ensuremath{\smash{\underset{\widetilde{\hphantom{#1}}}{#1}} \mathstrut}}

\newcommand{\domain}[1]{\ensuremath{\mathrm{dom}({#1})}}
\newcommand{\range}[1]{\ensuremath{\mathrm{codom}({#1})}}

\newcommand{\cal}{\ensuremath{\mathcal}}

\newtheorem {theorem1}{Theorem}[section]
\newtheorem {theorem}[theorem1]{Theorem}

\newtheorem {proposition}[theorem1]{Proposition}
\newtheorem {prop}[theorem1]{Proposition}
\newtheorem {lemma}[theorem1]{Lemma}
\theoremstyle{definition}
\newtheorem {definition}[theorem1]{Definition}
\theoremstyle{remark}
\newtheorem {remark}[theorem1]{Remark}

\newtheorem {example}[theorem1]{Example}

\newcommand{\calC}{\ensuremath{\mathcal{C}}}

\newcommand{\calO}{\ensuremath{\mathcal{O}}}

\newcommand{\calT}{\ensuremath{\mathcal{T}}}

\newcommand{\cat}[1]{\ensuremath{\EuScript #1}}

\newcommand{\Top}{\ensuremath{\calT\textup{op}}}

\newcommand{\ZZ}{\ensuremath{\mathbb{Z}}}

\DeclareMathOperator{\colim}{\textup{colim}}

\def\id{\ensuremath{\textsl{id}}}

\def\Res{\ensuremath{\textup{Res}}}

\newcommand{\rarrow}{\rightarrow}
\newcommand{\arrow}{\rightarrow}

\newcommand{\Poincare}{Poincar\'e~}

\theoremstyle{plain}

\begin{document}

\title[Equivariant cellular homology]{Equivariant cellular homology and its applications}
\author{Boris Chorny}
\address{Einstein Institute of Mathematics, Edmond Safra Campus, Givat Ram,
         The Hebrew University of Jerusalem, Jerusalem 91904, Israel }
\email{chorny@math.huji.ac.il}
\thanks{The author acknowledges the support of Sonderforschungsbereich (SFB) 478 of the
        University of M\"unster and Edmund Landau Center for Research in Mathematical Analysis.}

\subjclass{Primary 55N91; Secondary 55P91, 57S99}

\date{August 30, 2001}

\begin{abstract}
In this work we develop a cellular equivariant homology functor and apply
it to prove an equivariant Euler-\Poincare formula and an equivariant
Lefschetz theorem.
\end{abstract}

\maketitle

\section{Introduction}


  Let $D$ be an arbitrary small topologically enriched category.
  In this paper we develop a $D$-$CW$-homology functor which allows for
  easy computation of the \emph{ordinary} $D$-equivariant homology
  defined by E. Dror Farjoun in \cite{DF}.
  Our approach is a generalization of the $G$-$CW$-(co)homology functor
  constructed by S.J. Willson in \cite{Willson} for the case of $G$
  being compact Lie group.

\par  Then we apply the $D$-$CW$-homology functor to obtain:

\par \emph{(i)} Equivariant Euler-\Poincare formula:

  \begin{equation} \label{main}
     \chi^D(\dgrm X)=\sum_{n=0}^\infty(-1)^n \widetilde{rk}_\mathrm{HS}(H^D_n(\dgrm X; {\cal I}))
  \end{equation}

  This formula establishes a connection between the equivariant homology
  and an equivariant Euler characteristic; $\widetilde{rk}_\mathrm{HS}(\cdot)$
  is a slight modification of Hattori-Stallings rank (originally defined in \cite{Hatt},\cite{Stall}).

\par \emph{(ii)} Equivariant Lefschetz theorem:
    Let \dgrm X be a triangulated $D$-space, $f:\dgrm X \rightarrow
  \dgrm X$ an equivariant map.
  If the equivariant Lefschetz number
  \begin{equation}\label{main-lef}
     \Lambda_D (f) = \sum_{n=0}^\infty(-1)^n\widetilde{tr}_\mathrm{HS}(H^D_n(f;{\cal I}))
  \end{equation}
  is not equal to zero, then there are $f$-invariant orbits in \dgrm X,
  moreover the orbit types of the invariant orbits may be
  recovered from $\Lambda_D (f)$.


\paragraph{\emph{Acknowledgments.}}
\par    I would like to thank E.Dror Farjoun, V.Halperin, A.Libman, Sh. Rosset for  many
stimulating conversations and helpful ideas.

\section{Preliminaries}
\subsection{$D$-spaces.} Let \Top~ denote the category of the
  compactly generated Hausdorff topological spaces.
  Fix an arbitrary small category $D$ enriched over \Top.
  We work in the category $\Top^D$ of functors from $D$ to \Top.
  The objects of this category are called topological diagrams or just $D$-spaces.
  The arrows in $\Top^D$ are natural transformations
  of functors or \emph{equivariant maps}.

\subsection{$D$-homotopy.} An equivariant homotopy between two $D$-maps
$f,g:
  \dgrm X\arrow\dgrm Y$, where \dgrm X,\dgrm Y are $D$-diagrams, is a $D$-map
  $H:\dgrm X\times I\arrow\dgrm Y$, where $I$ denotes the constant $D$-space
  $I(d)=[0,1]$. A \emph{homotopy equivalence} $f:\dgrm X \arrow\dgrm Y$
  is a map with a (two sided) $D$-homotopy inverse.

\subsection{$D$-orbits.} We recall now the central concept of the $D$-homotopy
  theory (introduced in \cite{DF},\cite{DZ}) -- that of $D$-orbit.
  A $D$-orbit is a $D$-space $T:D\arrow\Top$, such that $\colim_D T=\{\ast\}$.
  A \textbf{free} $D$-orbit generated in $d\in ob(D)$ is $\Top^D \ni F^d=\hom_D (d,\cdot)$,
  i.e. $F^d(d')=\hom_D(d,d')$ and $F^d(d'\arrow d'')$ is given by
  the composition. Clearly $F^d$ is a $D$-orbit. A $D$-space \dgrm X is called free iff for any
  $s\in \colim_\cat C \dgrm X$ the full orbit $T_s$ lying over $s$ is free.

\subsection{$D$-$CW$-complexes.} A $D$-\emph{cell} is a $D$-space of the form
  $T\times e^n$, where T is a
  $D$-orbit and $e^n$ is the standard $n$-cell. An \emph{attaching map} of this
  $D$-cell to some $D$-space \dgrm X is a map $\phi:T\times\partial e^n\arrow\dgrm X$.
\par    A (\emph{relative}) $D$-$CW$-complex $(\dgrm X,\dgrm X_{-1})$ is a
  $D$-space \dgrm X together with a filtration
  $\dgrm X_{-1}\subset\dgrm X_0 \subset\ldots\subset \dgrm X_n \subset
  \dgrm X_{n+1}\subset\ldots\subset\dgrm X=\colim_n \dgrm X_n$,
  such that $\dgrm X_{n+1}$ is obtained from $\dgrm X_n$ by attaching a set of
  $n$-dimensional $D$-cells. Namely one has a push-out diagram of $D$-spaces:
\[\CD
        \coprod_i(T_i\times\partial e^n)  @>\phi>> \dgrm X_{n-1}\\
        @VVV                                     @VVV\\
        \coprod_i(T_i\times e^n)          @>\Phi>> \dgrm X_n
\endCD\]

\par  If $\dgrm X_{-1}=\varnothing \;$we call the $D$-$CW$-complex \emph{absolute}.

\par    Let \dgrm X be a $D$-$CW$-complex. A $D$-subspace $\dgrm Y \subset \dgrm X$
  is called the \emph{cellular subspace} if \dgrm Y has a $D$-$CW$-structure such
  that each cell of \dgrm Y is also a cell of \dgrm X.

\subsection{The category of orbits.}  The \emph{category of orbits} \cat O is a full
  topological sub-category of $\Top^D$ generated by all $D$-orbits.
\par    Usually \cat O is not a small category. For example for
  $D=J=(\bullet\!\!\longrightarrow\!\!\bullet)$,
  then $\cat O \cong \Top$. A model category has been constructed
  for the $D$-spaces of arbitrary orbit type in \cite{DF1}.
  We will be interested in the diagrams which are homotopy equivalent to
  the finite $D$-$CW$-complexes, i.e. only finite number of orbit types appear
  in such diagram. We collect those orbits into the a full
  subcategory $\calO'$ of $\cat O$ with a finite amount of objects.

\subsection{Orbit point $(\cdot)^\calO$ and realization $|\cdot|_D$ functors.}
        Suppose \calO~is a \emph{small} category of $D$-orbits. An orbit point functor
  $(\cdot)^\calO:\Top^D\arrow\Top^{\calO^{\mathrm{op}}}$ is a generalization
  to the diagram case of Bredon's fixed point functor. For any $D$-space \dgrm X
  (usually of type \calO) $(\dgrm X)^\calO$ is a $\calO^{\mathrm{op}}$ diagram
  such that $(\dgrm X)^\calO(T)=\hom_D(T,\dgrm X)$ for all $T\in ob(\calO)$ and
  the arrows of the diagram are induced by composition with the maps between
  orbits.
\par    If $f:\dgrm X\arrow\dgrm Y$ is an equivariant map between two
  $D$-spaces, then there exist an $\calO^\mathrm{op}$-equivariant map
  $f^\calO:\dgrm X^\calO\arrow\dgrm Y^\calO$, which is obtained from $f$ by
  composition:
\[
        \dgrm X^\calO(T)=\hom_D(T,\dgrm X) \ni g \overset{f^\calO}\longmapsto f\circ g \in
                         \hom_D(T,\dgrm Y)=\dgrm Y^\calO(T)
\]

\par    The fundamental property of $(\cdot)^\calO$ functor is that for any
  $D$-space \dgrm X the $\calO^\mathrm{op}$-space $(\dgrm X)^\calO$ is
  $\calO^\mathrm{op}$-free \cite[3.7]{DF}.

\begin{example}
        Consider a free $D$-space \dgrm X. And let the orbit category \calO~
  consists of all the free orbits. Then \calO~ is isomorphic to $D^\mathrm{op}$ as a category
  and $(\dgrm X)^\calO \cong \dgrm X$ (Yoneda's lemma).
\par    Another easy case occurs then \dgrm X is a $D$-$CW$-space. We
  shall discuss it in the next section.
\end{example}

\par    There exist a left adjoint to $(\cdot)^\calO$. It is called
  \emph{realization functor} $|\cdot|_D$, since it takes an
  $\calO^\mathrm{op}$-space and produce a $D$-space with the prescribed orbit
  point data (up to local weak equivalence). Realization functor in the group
  case has been constructed by A.D.Elmendorf in \cite{Elm} and has been
  generalized to the arbitrary diagram case by W.Dwyer and D.Kan in \cite{DK}.
  Compare also \cite{DF}.

\subsection{Equivariant Euler characteristic.}
    Let \dgrm X be a finite $D$-$CW$-complex, then \dgrm X is of
  type $\cal O$ for some category of orbits with finite amount of objects.
  We define the Equivariant Euler characteristic to be the
  Universal Additive Invariant \cite[I.5]{Lueck} $(U(D),\chi^D)$. Or,
  equivalently, we say that $\chi^D(\dgrm X)\in U(D)$ is equal to
  the alternating sum of the orbit types over the dimensions of the cells of \dgrm X
  in the free abelian group $U(D)=\bigoplus_{T\in Iso(h\calO)}\ZZ$ generated
  by the homotopy types of orbits.

\section{Equivariant cellular homology}

\subsection{$\calO^\mathrm{op}$-$CW$ structure on the orbit point space of
        a $D$-$CW$-complex.}
Let \calC~ denote a subcategory of $\Top^{\calO^\mathrm{op}}$
which is obtained as the image of $\Top^D$ under the functor
$(\cdot)^\calO$. Recall that there is the inclusion of  the
categories $\imath:D \hookrightarrow \calO^\mathrm{op}$, where
$\imath(d) = F^d$, for each $d\in ob(D)$. Hence there is a
functor $\Res:\Top^{\calO^\mathrm{op}}\rightarrow\Top^D$. By abuse
of notation we denote by $\Res$ also $\Res|_{\calC}$.
\begin{lemma}
The functor $(\cdot)^\calO$ is fully faithful.
\end{lemma}
\begin{proof}
The faithfulness is clear. We have to show only that for any map
$f_0:\dgrm X^\calO \rightarrow \dgrm Y^\calO$ there exists a map $f:\dgrm
X \rightarrow \dgrm Y$ s.t. $f_0 = f^\calO$. Take $f = \Res(f_0)$, then
if \dgrm X and \dgrm Y were orbits the result follows from the bijective
correspondence induced by the Yoneda's lemma: ${\mathrm hom}_D(\dgrm X,
\dgrm Y) = {\mathrm {hom}}_{\calO^{\mathrm {op}}}(\dgrm X^\calO, \dgrm
Y^\calO)$. The general claim will follow from the comparison of $f_0$ and
$f^\calO$ orbitwise, i.e. by their action on each full orbit. Fortunately
the functor $(\cdot)^\calO$, being right adjoint, commutes with taking
full orbit (pullback).
\end{proof}
\begin{lemma}
The pair of functors
$\Res(\cdot):\calC\leftrightarrow\Top^D:(\cdot)^\calO$ induce the
equivalence of the categories $\calC$ and $\Top^D$.
\end{lemma}
\begin{proof}
    We need to construct the natural isomorphisms of the functors
 $\id_{\Top^D} \cong \Res((\cdot)^\calO)$ and $\id_\calC\cong(Res(\cdot))^\calO$.
\par
    Let $\dgrm X \in \Top^D$, then $\Res(\dgrm X^\calO) \cong \dgrm X$
 because the generalized lemma of Yoneda \cite{Kelly} induces the objectwise homeomorphisms and
 the equivariance is preserved by the naturality of the Yoneda's isomorphism.
 But an equivariant map which is the objectwise homeomorphism is an isomorphism
 of $D$-spaces, hence the first isomorphism of functors.
\par
    Let $\dgrm X^\calO \in \calC$, then $\Res(\dgrm X^\calO) \cong \dgrm X$ by the
 first homeomorphism, then $(\Res(\dgrm X^\calO))^\calO \cong \dgrm X^\calO$.
 Hence the second isomorphism.{\sloppy

}
\end{proof}

\begin{proposition} \label{main:CW}
    Let \dgrm X be a (pointed) $D$-$CW$--space of orbit type \calO,where
  \calO is a small category of orbits. Consider the
  $\calO^\mathrm{op}$--space $\dgrm X^\calO$ to be the orbit point space
  of \dgrm X.
\par
    Then $\dgrm X^\calO$ has $\calO^\mathrm{op}$-$CW$ structure which
  corresponds to the $D$-$CW$ structure of \dgrm X in the following sense:
     let  $\dgrm pt_D\subseteq \dgrm X_0 \subseteq \dgrm X_1\subseteq \cdots
  \subseteq \dgrm X_n \subseteq \cdots \subseteq \dgrm X = {\colim_n} \dgrm X_n$
  is a $D$-$CW$ filtration of \dgrm X, such that each $\dgrm X_n$ is
  a push-out:
\[
    \CD
        {\coprod_i T_i\times S^{n-1}} @>{\phi}>> {\dgrm X_{n-1}}\\
            @VVV                 @VV{i_n}V\\
        {\coprod_i T_i\times D^{n}}   @>{\Phi}>>   {\dgrm X_{n}}
    \endCD
\]
then there exist a $\calO^\mathrm{op}$-$CW$--filtration:
  $\dgrm pt_{\calO^{\mathrm{op}}}\subseteq \dgrm X_0^\calO \subseteq \dgrm X_1^\calO
  \subseteq \cdots \subseteq \dgrm X_n^\calO \subseteq \cdots \subseteq
  \dgrm X^\calO= \colim_n \dgrm X_n^\calO$, such that
  $\dgrm X_n^\calO=(\dgrm X_n)^\calO$, and
\[
    \CD
        {\coprod_i F^{T_i}\times S^{n-1}} @>{\phi^\calO}>> {\dgrm X_{n-1}^\calO}\\
            @VVV                       @VV{i_n^\calO}V\\
        {\coprod_i F^{T_i}\times D^{n}}   @>{\Phi^\calO}>>   {\dgrm X_{n}^\calO}
    \endCD
\]
  is a push-out square.
\end{proposition}
\begin{proof}
    We proceed by the induction on the skeleton of \dgrm X.
\begin{eqnarray} \nonumber
    \dgrm X_0^\calO(T) = (\coprod T_i)^\calO(T) = \coprod((T_i)^\calO(T)) =
    \coprod\mathrm{hom}_{D}(T,T_i) = \\ \nonumber
    \coprod\mathrm{hom}_{\calO^\mathrm{op}}(F^{T},F^{T_i})= \coprod F^{T_i}(T).
\end{eqnarray}

    Hence the base of the induction.
\par
    Suppose we know the claim for $\dgrm X_n$. Then it follows
for $\dgrm X_{n+1}$ since $(\cdot)^\calO$ is both left and right adjoint,
so commutes both with push-outs and products.
\end{proof}

\subsection{$D$-$CW$-homology functor.}
  The construction of the (co)homology functor in \cite[4.16]{DF}
  depends on the specific $D$-$CW$-decomposition of $\dgrm X^{\cal O}$.
  We apply this construction to the cellular structure of $\dgrm X^\calO$,
  which was constructed in \ref{main:CW} and obtain the required
  $D$-$CW$-homology functor.
\subsection{Isotropy ring $\cal I$.}

In \cite{Willson} a universal coefficient system for the
$G$-equivariant homology have been developed(where $G$ is a
compact Lie group). Let us generalize this approach to the
coefficient systems for the classical $D$-homology theory.
Suppose $\cal O'$ is a small, full subcategory of the orbit
category $\cal O$. Let \dgrm X be a $D$-space of orbit type $\cal
O'$. Then a coefficient system for the ordinary (co)homology is a
homotopy (co)functor $M:{\cal O'}\rarrow (\Rmod)$.
\begin {definition}
Let R be a commutative ring. An \emph{isotropy ring } $\cal
I=I^{{\mathrm R,}O'}_{\mathrm D}$ is generated by ${\mathrm
mor}(h\cal O')$ as a free \Rmod. Define the multiplication on the
generators by
\[
fg =\left\{
    \begin{array}{rl}
        f \mathbin{\circ} g,  & \textup {if\;} \range g = \domain f \\
        0,                    & \textup {otherwise}
    \end{array}
    \right.
\]
and extend the definition to the rest of the elements of $\cal I$
by linearity.
\end{definition}
\begin{prop}
The category $\cal M$ of the left $\cal I$-modules which satisfy:
\begin{equation} \label{cond}
    \forall M\in ob({\cal M}),\; M=\bigoplus_{T\in ob(h{\cal O}')}{1_T M}
\end{equation}
(where $\{1_T M\}_{T\in ob(h{\cal O}')}$ are left $R$-modules) and
the category of $R(h{\cal O'}){\textup -mod}$ of functors from
$h{\cal O'}$ to the category of left $R$-modules are equivalent.
\end{prop}
\begin{proof}
Let us define a pair of functors which induce the required equivalence:
\[
    \zeta: {\cal M} \leftrightarrows R(h{\cal O'}){\textup -mod}:\xi.
\]
Let $M\in ob({\cal M})$, $T \in ob(h{\cal O'})$, then define
\[
    \zeta M(T) = 1_T M.
\]
If ${\textup{mor}}(h{\cal O'})\ni f:T_1\rightarrow T_2$, then
define
\[
    \zeta M(f)(1_{T_1}m) = f1_{T_1}m = (1_{T_2}f)1_{T_1}m \in
    1_{T_2}M.
\]
Obviously the morphisms of the left $\cal I$ modules correspond
to the natural transformations of the functors.
\par
    Given a $R(h{\cal O'})$-module $N$, then
\[
    \xi N = \bigoplus_{T\in ob(h{\cal O'})}N(T), \textup{\;as\;
    a\; left\;} R\textup{-module.}
\]
Define the left $\cal I$-module structure on $\xi N$ by
$f(\ldots,n,\dots) = (\ldots,fn,\dots)$, where
\[
   N(\range f)\ni
   fn = \left\{
            \begin{array}{rl}
            f(n),   & \textup{if\;} n \in N(\domain f)\\
            0,      & \textup{otherwise}
            \end{array}
            \right.
\]
Now it is clear that the defined functors provide the equivalence of the
categories.
\end{proof}

\begin{remark}\label{rem1}
The ring $\cal I$ considered as a left $\cal I$-module is an
object of $\cal M$, because ${\cal I} \cong \bigoplus_{T\in
ob(h{\cal O}')}{1_T {\cal I}}$ (as left $R$-modules) by the
construction. But it also carries an obvious structure of the
right $\cal I$ module, so the $\zeta {\cal I}(T)$.
\end{remark}
\begin{remark}\label{rem2}
If $ob(h\cal O')$ is finite then the ring $\cal I$ has a
two-sided identity element $1 = \sum_{T\in ob(h{\cal O'})}1_T$
together with its decomposition into the sum of the orthogonal
idempotents and the condition (\ref{cond}) is redundant.
\end{remark}
\begin{definition}
The augmentation $\varphi:{\cal I} \rightarrow \bigoplus_{T\in
Iso(ob(h{\cal O'}))}R$ is defined for any
\[
    {\cal I} \ni g = \sum_{T\in ob(h{\cal O'})}\sum_{f \in \mathrm{mor}(T,T)} r_{f} f + \sum_{h \in \mathrm{mor}(T_1,T_2),T_1 \neq T_2} s_h h
\]
(only a finite number of $r_{f},s_{h} \in R$ is non equal to
zero) to be
\[
    \varphi(g)=(\ldots, \sum_{f \in \mathrm{mor}(T,T)} r_{f} ,\ldots)\in\bigoplus_{T\in Iso(ob(h{\cal O'}))}R
\]
\end{definition}
\begin{remark} \label{rem3}
The idempotents in $\cal I$ which correspond to the $D$-homotopy
equivalent orbits are identified under $\varphi$. Apparently,
$\varphi$ is an epimorphism of rings. Consider the abelinization
functor $Ab:(Rings) \rightarrow \cat{Ab}$ which corresponds to a
ring its additive group divided by the commutator subgroup. Then
$Ab(\varphi):Ab({\cal I})\rightarrow \bigoplus_{T\in
Iso(ob(h{\cal O'}))}R$. The last map will be used to obtain a
generalization of the Euler-\Poincare formula.
\end{remark}

\section{Applications}
    Let \dgrm X be a finite $D$-$CW$-complex of type $\cal O'\;$
for some orbit category $\cal O'$ with $\mathrm{obj}(\cal O')$ a
\emph{finite} set.

\subsection{Equivariant Euler-\Poincare formula.}
\par    We remind that the equivariant Euler characteristic lies in the
abelian group $U(D)\cong \bigoplus_{Iso(obj(h{\cal O'}))}\ZZ$, so
in order to apply Hattori--Stallings machinery we need to choose a
coefficient system for the equivariant homology such that the
resulting chain complex and homology groups will be endowed with
the module structure over some ring $S$ which allows an
epimorphism $\varepsilon: Ab(S)\longrightarrow U(D)$.

\par Our choice of the coefficient system for the equivariant
homology will be the isotropy ring $\cal I = I^{\ZZ,O'}_{\mathrm D
}$ taken over itself as a left module.

\begin{lemma} \label{main:lemma}
Let \dgrm X be a finite $D$-$CW$-complex. Suppose \dgrm X has
$n_q$ $q$-dimensional cells and $t_1+\cdots+t_s = n_q$, $t_i$ is
the number of $q$-dimensional cells of the same homotopy type $T_i
\in Iso(obj(h{\cal O'}))$. Then $ \calC_q(\dgrm X) \otimes_{\cal
O'} {\zeta{\cal I}} \cong \zeta{\cal I}(T_1) ^{t_1} \oplus \cdots
\oplus \zeta{\cal I}(T_n)^{t_s}$ as a left \ZZ-module.
\end{lemma}
\begin{proof}
Let $t_i = r_{i1}+\cdots+r_{ik}$, where $r_{ij}$ is the number of
$q$-dimensional cells of type $T_{ij}\in ob(\calO)$ of homotopy type
$T_i$. By the construction of the equivariant homology $\calC_q(\dgrm X)
= \bigoplus_{i=1}^s( \oplus_{j=1}^k \ZZ(\mathrm{hom}_{\cal O'}(?,
T_{ij})^{r_{ij}}))$. The dual Yoneda isomorphism \cite[p.74]{Kelly}
implies:
\[
    \calC_q(\dgrm X)\otimes_{\calO'} \zeta {\cal I}\cong \bigoplus_{i=1}^s( \oplus_{j=1}^k
    \zeta{\cal I}( T_{ij})^{r_{ij}}) \cong \bigoplus_{i=1}^s\oplus_{j=1}^k(1_{T_{ij}}{\cal I})^{r_{ij}} ,
\]
If $T_{ij_1}$ is isomorphic to $T_{ij_2}$ in $h\cal O'$ then
there is an obvious isomorphism of the left $\ZZ$-modules and
right $\cal I$-modules $1_{T_{ij_1}}{\cal I} \cong 1_{T_{ij_2}}
{\cal I}$. Let us choose a representative $T_i$ of each
isomorphism class of objects in $h\cal O'$, then
\[
\calC_q(\dgrm X)\otimes_{\calO'} \zeta {\cal I} \cong
\bigoplus_{i=1}^s(1_{T_{i}}{\cal I})^{(\sum_{j=1}^k r_{ij})} \cong
\bigoplus_{i=1}^s(1_{T_{i}}{\cal I})^{t_{i}}\cong
\bigoplus_{i=1}^s (\zeta{\cal I}( T_{i}))^{t_{i}}
\]
\end{proof}
Because of \ref{rem1} the equivariant chain complex $\{\calC_q (\dgrm X)
\otimes_{\calO'} \zeta {\cal I}\}_{q = 0}^{\mathrm{dim} \dgrm X}$ is a
complex of \emph{projective} right $\cal I$-modules and the equivariant
homology is endowed with the right $\cal I$-module structure.
\paragraph{\emph{Notation:}} $\chi_{_\mathrm{HS}}(\cdot)$ means Euler
characteristic of a $\cal I$ differential complex with respect to
$rk_\mathrm{HS}(\cdot)$.
\begin{prop} \label{intermidiate:prop}
Let $K_\ast=\calC_\ast (\dgrm X) \otimes_{\cal O'} \zeta{\cal I}$
be a right $\cal I$-complex, then $\chi^D(\dgrm X) =
Ab(\varphi)(\chi_{_\mathrm{HS}}(K_\ast))$ whenever left side is
defined.{\sloppy

}
\end{prop}
\begin{proof}
It is easy to see that $rk_{\mathrm HS}(1_T{\cal I}) = 1_T \in Ab(\cal
I)$. Lemma \ref{main:lemma} together with \ref{rem3} completes the proof.
\end{proof}
Now we combine \ref{intermidiate:prop} with the additivity properties of
the Hattory-Stallings rank and obtain the following
\begin{theorem}
$\chi^D(\dgrm X) = Ab(\varphi)(\sum_{n=0}^\infty (-1)^n
rk_\mathrm{HS} H_n^D(\dgrm X;\zeta {\cal I}))$, whenever the left
side is defined.
\end{theorem}
\begin{example}
Consider the $J$-diagram:
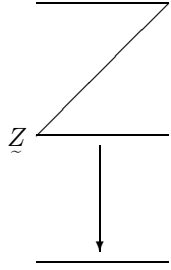
\begin{center}
\begin{picture}(51,100)
    \put(-10,45){\dgrm Z}

    \put(1,99){\line(1,0){50}}     
    \put(51,99){\line(-1,-1){50}}  
    \put(1,49){\line(1,0){50}}     

    \put(25,45){\vector(0,-1){40}}  

    \put(1,1){\line(1,0){50}}      
\end{picture}
\end{center}

\dgrm Z has two $0$-cells of type
$T_2=[\overset{..}{\underset{\cdot}{\downarrow}}]$ and one $1$-cell of type
$T_3=[\overset{...}{\underset{\cdot}{\downarrow}}]$, hence \mbox{$ \chi^J(\dgrm
Z)=2[\underset{\cdot}{\overset{..}{\downarrow}}] -
[\underset{\cdot}{\overset{...}{\downarrow}}]$}.{\sloppy

}

\par    The category $\cal O'$ of orbits contains two objects: $T_2$, $T_3$.
The cellular chain complex tensored with the coefficients ${\cal I
} = {\cal I}_J^{\ZZ,\{T_2,T_3\}}$ becomes:
\[
    \cdots\arrow 0\arrow 1_{T_3}{\cal I} \overset{\partial_1}{\arrow} (1_{T_2}{\cal I})^2
\]
and $\partial_1=0$ from the orbit type considerations.
\par    $U(J)=\ZZ\oplus\ZZ$ in that case. And $H_0^J(\dgrm Z,{\cal I})=(1_{T_2}{\cal I})^2$,
$H_1^J(\dgrm Z,{\cal I})=1_{T_3}{\cal I}$ are right $\cal I$ - modules.
Hence, $\chi^J(\dgrm Z)=(2,0)-(0,1)=(2,-1)$.
\par    Let us, for comparison, calculate the $J$-equivariant homology of
  \dgrm Z with $\ZZ^{\cal O'}$ coefficients: $H_i^J(\dgrm Z,\ZZ^{\cal O'})=
  H_i(\colim_J \dgrm Z,\ZZ)$ (see \cite[5.2]{DF}). Then $\colim_J \dgrm Z=I=
  [0,1]$ and
\[
    H_i^J(\dgrm Z,\ZZ^{\cal O'})=\left\{
                    \begin{array}{ll}
                    \ZZ,    &i=0\\
                    0,  &\text{otherwise}
                    \end{array}
                \right.
\]
We can see that $\ZZ^{\cal O'}$ coefficients are inappropriate to the
Euler-\Poincare formula.
\end{example}

\subsection {Equivariant Lefschetz theorem.}
\par
        Using cellular equivariant homology functor we are
able now to proof a version of the equivariant Lefschetz
theorem.
\par
    Some result of the Lefschetz type in the equivariant setting may
be obtained already by applying the ordinary Lefschetz theorem:
consider an equivariant map $f: \dgrm X \rightarrow \dgrm X$,
where \dgrm X is a diagram over small category $D$, then if the
Lefschetz number $\Lambda (\colim_D \dgrm X)\neq 0$ there are $f$-
invariant $D$-orbits in \dgrm X. However the advantage of using
the equivariant homology and equivariant Lefschetz number
$\Lambda_D(\dgrm X) \in U(D)$ is that we obtain the specific
information about orbit type of the invariant orbit.

\par First we give a technical
\begin {definition}
        A $D$-$CW$-complex \dgrm X will be called the triangulated $D$-space if the natural
        $CW$-structure of $\colim {\dgrm X}$ also triangulates $\colim {\dgrm X}$.
\end {definition}

  The following lemma will be used in the proof of the equivariant Lefschetz theorem.

\begin {lemma} \label {refinement:lemma}
        Let \dgrm X be a triangulated diagram, then for any refinement $Y$ of the triangulation of
        $\colim \dgrm X$, there exists a $D$-$CW$-complex $\dgrm X'$, such that $\dgrm X'$ is
        $D$-homeomorphic to \dgrm X and $\colim \dgrm X' = Y$ (as the triangulated spaces).
        $\dgrm X'$ will be called the refinement of \dgrm X.
\end {lemma}

\begin{proof}
Consider a new simplex $\Delta$ in the triangulation of $Y$. It lies in
some old simplex of $\colim \dgrm X :\Delta \in \Delta'$. Then consider
the pull-back:
\[
        \lim \left(
                \begin{array}{ccc}
                        &                       & \dgrm X    \\
                        &                       & \downarrow \\
                \Delta  &\hookrightarrow        & \colim_D \dgrm X
                \end{array}  \right)
            = T\times\Delta,
\]
  where $T$ is the orbit which lies over $\Delta'$.
\par
We've obtained the cell of the new $D$-$CW-$complex $\dgrm X'$.
Continuing in the same way for the rest of the simplices of $Y$ completes
the construction of $\dgrm X'$. Hence $D$-$CW$-complex $\dgrm X'$ has the
same underlying topological diagram as $\dgrm X$, therefore they are
$D$-homeomorphic.
\end{proof}

\begin {definition}
Let $f:\dgrm X \rightarrow \dgrm X$ be a map of the finite
triangulated $D$-space \dgrm X of orbit type $\cal O'$, where
$\cal O'$ is an orbit category with the finite number $n$ of
objects. Let ${\cal I} = {\cal I}_D^{\ZZ, {\cal O}'}$. Then the
equivariant Lefschetz number of $f$:
\[
U(D)\owns(\lambda_1,\ldots,\lambda_n) = \Lambda_D(f) =
Ab(\varphi)(\sum_{k=0}^\infty(-1)^k\traceHS(H_k(f;{\cal I})))
\]
\end {definition}

\begin {theorem}
Let \dgrm X be a finite triangulated diagram over $D$. $f:\dgrm X
\rightarrow \dgrm X$ be a $D$-map. $\Lambda_D(f) = (\lambda_1, \ldots,
\lambda_n)\in U(D)$ - Lefschetz number of $f$. Then if there is no
$f$-invariant orbit of type $T_m$, $\lambda_m=0$.
\end {theorem}
\begin{proof}
A simplex in $\colim \dgrm X$ will be called of type $T$ if the
overlaying orbit is of type $T$ in \dgrm X. Then the condition that there
are no invariant orbits of type $T_m$ is equivalent to the condition that
there are no fixed points in the simplices of type $T_m$.
\par
Since \dgrm X is a finite triangulated diagram, $\colim \dgrm X$ is a
finite triangulated space, hence it is a compact metric space. If there
are no fixed points of type $T_m$, then there exists a refinement $Y$ of
the triangulation such that if $\Delta$ is a simplex of type $T_m$ in
$Y$, $\Delta \cap (\colim f)(\Delta) = \varnothing$.
\par
Consider the refinement $\dgrm X'$ of \dgrm X, which exists by lemma
\ref{refinement:lemma}. Since $\dgrm X' \cong \dgrm X$, $H^D_\ast(\dgrm
X';{\cal I}) = H^D_\ast(\dgrm X;{\cal I})$,
 $\Lambda(f') = \Lambda(f)$, where $f':\dgrm X'\rightarrow \dgrm X'$ is
equal to $f$, $\lambda_m' = \lambda_m$. Therefore, it is enough
to show that $\lambda_m'=0$.{\sloppy

}
\par
Now,
\[
\Lambda_D(f') =
Ab(\varphi)(\sum_{k=0}^\infty(-1)^k\traceHS(H_k(f';{\cal I}))) =
Ab(\varphi)(\sum_{k=0}^\infty(-1)^k\traceHS({\cal C}_k(f';{\cal
I}))),
\]
where ${\cal C}_k(f';{\cal I})$ is the map induced by $f$ on the
chains ${\cal C}_k(\dgrm X;{\cal I}) =
 {\cal C}_k(\dgrm X)\otimes_{\cal O'}{\cal I} =
 (1_{T_1}{\cal I})^{t_1}\oplus\ldots\oplus (1_{T_n}{\cal I})^{t_n}$
as $\cal I$-module. Because of the property: $\Delta\cap(\colim
f)(\Delta) = \varnothing$ for \emph{any} simplex $\Delta$ of type $T_m$,
the induced map on ${\cal C}_k(\dgrm X;{\cal I})$ will take the generator
$1_{T_m}$ corresponding to $\Delta$ outside the submodule $1_{T_m}{\cal
I}$, that it generates. Then the $m$-th entry of
$Ab(\varphi)(\sum_{k=0}^\infty(-1)^k\traceHS({\cal C}_k(f';{\cal I})))$
will be zero. This is true for all $k$, hence $\lambda_m = 0$.
\end{proof}

\end{document}